\documentclass[onecolumn]{autart}
\usepackage{eufrak,mathrsfs,amsmath,longtable}
\input{amssym}
\newcommand{\g}{\goth g}

\newcommand{\p}{\rm \partial}

\newcommand{\EE}{\mathscr E}

\renewcommand{\theequation}{\arabic{section}.\arabic{equation}}
\begin{document}
\begin{frontmatter}

\title{Preliminarily  group classification  of  a class of \\ 2D nonlinear
heat equations
} 
\thanks[footnoteinfo]{ Corresponding author: Tel. +9821-73913426.
Fax +9821-77240472.}
\author[]{M. Nadjafikhah\thanksref{footnoteinfo}}\ead{m\_nadjafikhah@iust.ac.ir},
\author[]{R. Bakhshandeh-Chamazkoti}\ead{r\_bakhshandeh@iust.ac.ir},
\address{School of Mathematics, Iran University of Science and Technology, Narmak, Tehran 1684613114, Iran.}
\renewcommand{\sectionmark}[1]{}
\begin{abstract}
A preliminary group classification of the class 2D nonlinear heat
equations $u_t=f(x,y,u,u_x,u_y)(u_{xx}+u_{yy})$, where $f$  is
arbitrary smooth function of the variables $x,y,u,u_x$ and $u_y$
using  Lie method, is given. The paper is one of the few
applications of an algebraic approach to the problem of group
classification: the method of preliminary group classification.
\end{abstract}
\end{frontmatter}
\section{Introduction}
It is well known that the symmetry group method plays an important
role in the analysis of differential equations.   The history of
group classification methods goes back to Sophus Lie. The first
paper on this subject is \cite{[1]}, where Lie proves that a
linear two-dimensional second-order PDE may admit at most a
three-parameter invariance group (apart from the trivial infinite
parameter symmetry group, which is due to linearity).  He computed
the maximal invariance group of the one-dimensional heat
conductivity equation and utilized this symmetry to construct its
explicit solutions. Saying it the modern way, he performed
symmetry reduction of the heat equation. Nowadays symmetry
reduction is one of the most powerful tools for solving nonlinear
partial differential equations (PDEs). Recently, there have been
several generalizations of the classical Lie group method for
symmetry reductions. Ovsiannikov \cite{[2]} developed the method
of partially invariant solutions. His approach is based on the
concept of an equivalence group, which is a Lie transformation
group acting in the extended space of independent variables,
functions and their derivatives, and preserving the class of
partial differential equations under study. In an attempt to study
nonlinear effects Saied and Hussain \cite{[3]} gave some new
similarity solutions of the (1+1)-nonlinear heat equation. Later
Clarkson and Mansfield \cite{[4]} studied classical and
nonclassical symmetries of the (1+1)-heat equation and gave new
reductions for the linear heat equation and a catalogue of
closed-form solutions for a special choice of the function
$f(x,y,u,u_x,u_y)$ that appears in their model. In higher
dimensions Servo \cite{[5]} gave some conditional symmetries for a
nonlinear heat equation while Goard et al. \cite{[6]} studied the
nonlinear heat equation in the degenerate case. Nonlinear heat
equations in one or higher dimensions are also studied in
literature by using both symmetry as well as other methods
\cite{[7],[8]}.\newline
There are a number of papers
to study (1+1)-nonlinear heat equations from the point of view of
Lie symmetries method. The (2+1)-dimensional nonlinear heat
equations
\begin{eqnarray}
u_t=f(u)(u_{xx}+u_{yy}),\label{eq:1}
\end{eqnarray}
are investigated in \cite{[9]} and in present paper we studied
\begin{eqnarray}
u_t=f(x,y,u,u_x,u_y)(u_{xx}+u_{yy}),\label{eq:2}
\end{eqnarray}
Similarity techniques are applied in \cite{[10],[11],[12],[13]} for
(2+1)-dimensional wave equations.
\section{Symmetry Methods}
Let a partial differential equation contains  $p$ dependent
variables and $q$ independent variables. The one-parameter Lie
group of transformations
\begin{eqnarray}
x_i\longmapsto x_i+\epsilon\xi^i(x,u)+O(\epsilon^2);\hspace{1cm}
u_{\alpha}\longmapsto
u_{\alpha}+\epsilon\varphi^{\alpha}(x,u)+O(\epsilon^2),\label{eq:3}
\end{eqnarray}
where $i=1,\ldots,p$ and $\alpha=1,\ldots,q$. The action of the
Lie group can be recovered from that of its associated
infinitesimal generators. we consider general vector field
\begin{eqnarray}
X=\sum_{i=1}^p\xi^i(x,u)\frac{\p}{{\p}x_i}+
\sum_{\alpha=1}^q\varphi^{\alpha}(x,u)\frac{\p}{{\p}u^{\alpha}}.\label{eq:4}
\end{eqnarray}
on the space of independent and dependent variables.
%
%
%
The symmetry generator associated with (\ref{eq:4}) given by
\begin{eqnarray}
X=\xi^1(x,y,t,u)\frac{\p}{{\p}x}+\xi^2(x,y,t,u)\frac{\p}{{\p}y}+
\xi^3(x,y,t,u)\frac{\p}{{\p}t}+\varphi(x,y,t,u)\frac{\p}{{\p}u}.\label{eq:6}
\end{eqnarray}
The second prolongation of $X$ is the vector field
\begin{eqnarray}\nonumber
 X^{(2)}=X+\varphi^x\frac{\p}{{\p}u_x}+\varphi^y\frac{\p}{{\p}u_y}+\varphi^t\frac{\p}{{\p}u_t}+
\varphi^{xx}\frac{\p}{{\p}u_{xx}}+\varphi^{xy}\frac{\p}{{\p}u_{xy}}
+\varphi^{xt}\frac{\p}{{\p}u_{xt}}
+\varphi^{yy}\frac{\p}{{\p}u_{yy}}
+\varphi^{yt}\frac{\p}{{\p}u_{yt}}+\varphi^{tt}\frac{\p}{{\p}u_{tt}},\\\label{eq:7}
\end{eqnarray}
that its coefficients are obtained with following formulas
\begin{eqnarray}\label{eq:8}
&&\varphi^x={D}_x\varphi-u_x{D}_x\xi^1-u_y{D}_x\xi^2-u_t{D}_x\xi^3,\hspace{2cm}
\varphi^y={D}_y\varphi-u_x{D}_y\xi^1-u_y{D}_y\xi^2-u_t{D}_y\xi^3,\\\nonumber
&&\varphi^t={D}_t\varphi-u_x{D}_t\xi^1-u_y{D}_t\xi^2-u_t{\rm
D}_t\xi^3,\hspace{2.1cm}
\varphi^{xx}={D}_x\varphi^x-u_{xx}{D}_x\xi^1-u_{xy}{D}_x\xi^2-u_{xt}{D}_x\xi^3\\\nonumber
&&\hspace{-2mm}\varphi^{yy}={D}_y\varphi^y-u_{xy}{D}_y\xi^1-u_{yy}{D}_y\xi^2-u_{yt}{D}_y\xi^3\hspace{1.58cm}
\varphi^{tt}={D}_t\varphi^t-u_{xt}{D}_t\xi^1-u_{yt}{D}_t\xi^2-u_{tt}{D}_t\xi^3\\\nonumber
&&\hspace{-2mm}\varphi^{yt}={D}_y\varphi^t-u_{xy}{D}_y\xi^1-u_{yy}{D}_y\xi^2-u_{yt}{D}_y\xi^3\hspace{1.58cm}
\varphi^{xt}={D}_t\varphi^t-u_{xt}{D}_t\xi^1-u_{yt}{D}_t\xi^2-u_{tt}{D}_t\xi^3
\end{eqnarray}
where the operators $D_x$, $D_y$ and $D_t$ denote the total
derivatives with respect to $x,y$ and $t$:
\begin{eqnarray}\nonumber
D_x&=&\frac{\p}{{\p}x}+u_x\frac{\p}{{\p}u}+u_{xx}\frac{\p}{{\p}u_x}+u_{xy}\frac{\p}{{\p}u_y}+
u_{xt}\frac{\p}{{\p}u_t}+\ldots\\
D_y&=&\frac{\p}{{\p}y}+u_y\frac{\p}{{\p}u}+u_{yy}\frac{\p}{{\p}u_y}+u_{yx}\frac{\p}{{\p}u_x}+
u_{yt}\frac{\p}{{\p}u_t}+\ldots\\\label{eq:9}
D_t&=&\frac{\p}{{\p}t}+u_t\frac{\p}{{\p}u}+u_{tt}\frac{\p}{{\p}u_t}+u_{tx}\frac{\p}{{\p}u_x}+
u_{ty}\frac{\p}{{\p}u_y}+\ldots\nonumber
\end{eqnarray}
%
%
%
By theorem 6.5. in \cite{[14]},
$X^{(2)}[u_t-f(x,y,u,u_x,u_y)(u_{xx}+u_{yy})]|_{(6)}=0$ whenever
\begin{eqnarray}
u_t-f(x,y,u,u_x,u_y)(u_{xx}+u_{yy})=0.\label{eq:10}
\end{eqnarray}
Since
$$X^{(2)}[u_t-f(x,y,u,u_x,u_y)(u_{xx}+u_{yy})]=\varphi^t-
(f_x\xi^1+f_y\xi^2+f_u\varphi+f_{u_x}\varphi^x+f_{u_y}\varphi^y)
(u_{xx}+u_{yy})-f(x,y,u,u_x,u_y)(\varphi^{xx}+\varphi^{yy}),$$
therefore we obtain the following determining function:
\begin{eqnarray}
\varphi^t-(f_x\xi^1+f_y\xi^2+f_u\varphi+f_{u_x}\varphi^x+f_{u_y}\varphi^y)
(u_{xx}+u_{yy})
-f(x,y,u,u_x,u_y)(\varphi^{xx}+\varphi^{yy})=0.\label{eq:11}
\end{eqnarray}
In the case of arbitrary $f$ it follows
\begin{eqnarray}
\xi^1=\xi^2=\varphi=0,\label{eq:12}
\end{eqnarray}
or
\begin{eqnarray}
\xi^1=\xi^2=\varphi=0,\;\;\;\;\;\xi^3=C.\label{eq:13}
\end{eqnarray}
Therefore, for arbitrary $f(x,y,u,u_x,u_y)$ Eq. (\ref{eq:1})
admits the one-dimensional Lie algebra ${\g}_1$, with the basis
\begin{eqnarray}
X_1=\frac{\p}{{\p}t}.\label{eq:14}
\end{eqnarray}
${\g}_1$ is called the principle Lie algebra for Eq. (\ref{eq:1}).
So, the remaining part of the group classification is to specify
the coefficient $f$ such that Eq. (\ref{eq:1}) admits an extension
of the principal algebra ${\g}_1$. Usually, the group
classification is obtained by inspecting the determining equation.
But in our case the complete solution of the determining equation
(\ref{eq:11}) is a wasteful venture. Therefore, we don't solve the
determining equation but, instead we obtain a partial group
classification of Eq. (\ref{eq:1}) via the so-called method of
preliminary group classification. This method was suggested in
\cite{[10]} and applied when an equivalence group is generated by
a finite-dimensional Lie algebra ${\g}_{\EE}$. The essential part
of the method is the classification of all nonsimilar subalgebras
of ${\g}_{\EE}$. Actually, the application of the method is simple
and effective when the classification is based on
finite-dimensional equivalence algebra ${\g}_{\EE}$.

\section{Equivalence transformations}
An equivalence transformation is a nondegenerate change of the
variables $t,x,y,u$ taking any equation of the form (\ref{eq:1})
into an equation of the same form, generally speaking, with
different $f(x,y,u,u_x,u_y)$. The set of all equivalence
transformations forms an equivalence group ${\EE}$. We shall find
a continuous subgroup ${\EE}_C$ of it making use of the
infinitesimal method.

We consider an operator of the group ${\EE}_C$ in the form
\begin{eqnarray}
Y=\xi^1(x,y,t,u)\frac{\p}{{\p}x}+\xi^2(x,y,t,u)\frac{\p}{{\p}y}+
\xi^3(x,y,t,u)\frac{\p}{{\p}t}+\varphi(x,y,t,u)\frac{\p}{{\p}u}
+\mu(x,y,t,u,u_x,u_y,u_t,f)\frac{\p}{{\p}f},\label{eq:15}
\end{eqnarray}
from the invariance conditions of Eq. (\ref{eq:1}) written as the
system:
\begin{eqnarray}\label{eq:3-16}
u_t&-&f(x,y,u,u_x,u_y)(u_{xx}+u_{yy})=0,\\\nonumber
f_t&=&f_{u_t}=0,
\end{eqnarray}
where $u$ and $f$ are considered as differential variables: $u$ on
the space $(x,y,t)$ and $f$ on the extended space
$(x,y,t,u,u_x,u_y)$.

The invariance conditions of the system (\ref{eq:3-16})  are
\begin{eqnarray}\label{eq:17}
Y^{(2)}(u_t&-&f(x,y,u,u_x,u_y)(u_{xx}+u_{yy}))=0,\\\nonumber
Y^{(2)}(f_t)&=&Y^{(2)}(f_{u_t})=0,
\end{eqnarray}
where $Y^{(2)}$ is the prolongation of the operator (\ref{eq:15}):
\begin{eqnarray}\label{eq:3-18}
Y^{(2)}=Y+\varphi^x\frac{\p}{{\p}u_x}+\varphi^y\frac{\p}{{\p}u_y}+\varphi^t\frac{\p}{{\p}u_t}+
\varphi^{xx}\frac{\p}{{\p}u_{xx}}+\varphi^{xy}\frac{\p}{{\p}u_{xy}}
+\varphi^{xt}\frac{\p}{{\p}u_{xt}}&+&\varphi^{yy}\frac{\p}{{\p}u_{yy}}\\\nonumber
&+&\varphi^{yt}\frac{\p}{{\p}u_{yt}}+\varphi^{tt}\frac{\p}{{\p}u_{tt}}+
\mu^t\frac{\p}{{\p}f_{t}}+\mu^{u_t}\frac{\p}{{\p}f_{u_t}}.
\end{eqnarray}
The coefficients $\varphi^x, \varphi^y, \varphi^t, \varphi^{xx},
\varphi^{xy}, \varphi^{xt}, \varphi^{yy}, \varphi^{yt},
\varphi^{tt}$ are given in (\ref{eq:8}) and the other coefficients
of (\ref{eq:3-18}) are obtained by applying the prolongation
procedure to differential variables $f$ with independent variables
$(x,y,t,u,u_x,u_y,u_t)$. we have
\begin{eqnarray}
\mu^t&=&\widetilde{D}_t(\mu)-f_x\widetilde{D}_t(\xi^1)-f_y\widetilde{D}_t(\xi^2)
-f_u\widetilde{D}_t(\varphi)-f_{u_x}\widetilde{D}_t(\varphi^x)-f_{u_y}\widetilde{D}_t(\varphi^y),\label{eq:19}\\
\mu^{u_t}&=&\widetilde{D}_{u_t}(\mu)-f_x\widetilde{D}_{u_t}(\xi^1)-f_y\widetilde{D}_{u_t}(\xi^2)
-f_u\widetilde{D}_{u_t}(\varphi)-f_{u_x}\widetilde{D}_{u_t}(\varphi^x)-f_{u_y}\widetilde{D}_{u_t}(\varphi^y),\label{eq:20}
\end{eqnarray}
where
\begin{eqnarray}
\widetilde{D}_t=\frac{\p}{{\p}t},\hspace{1cm}\widetilde{D}_{u_t}=\frac{\p}{{\p}u_t}.\label{eq:21}
\end{eqnarray}
So, we have the following prolongation formulas:
\begin{eqnarray}\label{eq:22}
\mu^t&=&\mu_t-f_x\xi_t^1-f_y\xi_t^2-f_u\varphi_t-f_{u_x}(\varphi^x)_t-f_{u_y}(\varphi^y)_t,\\\nonumber
\mu^{u_t}&=&\mu_{u_t}-f_{u_x}(\varphi^x)_{u_t}-f_{u_y}(\varphi^y)_{u_t},
\end{eqnarray}
By the invariance conditions (\ref{eq:17}) give rise to
\begin{eqnarray}
\mu^t=\mu^{u_t}=0,\label{eq:23}
\end{eqnarray}
that is hold for every $f$. Substituting (\ref{eq:23}) into
(\ref{eq:22}), we obtain
\begin{eqnarray}
\begin{array}{ll}
\mu_t=\mu_{u_t}=0\\
\xi^1_x=\xi^2_t=\varphi_t=0\\
(\varphi^x)_t=(\varphi^x)_{u_t}=(\varphi^y)_t=(\varphi^y)_{u_t}=0\label{eq:24}
\end{array}
\end{eqnarray}
Moreover with substituting (\ref{eq:3-18}) into (\ref{eq:17}) we
obtain
\begin{eqnarray}
\varphi^t-f(x,y,u,u_x,u_y)(\varphi^{xx}+\varphi^{yy})-\mu(u_{xx}+u_{yy})=0.\label{eq:25}
\end{eqnarray}
We are left with a polynomial equation involving the various
derivatives of $u(x,y,t)$ whose coefficients are certain
derivatives of $\xi^1,\xi^2,\xi^3$ and $\varphi$. Since
$\xi^1,\xi^2,\xi^3,\varphi$ only depend on $x,y,t,u$ we can equate
the individual coefficients to zero, leading to the complete set
of determining equations:
\begin{eqnarray}
\xi^1&=&\xi^1(x,y)\label{eq:26}\\
\xi^2&=&\xi^2(y)\label{eq:27}\\
\xi^3&=&\xi^3(t)\label{eq:28}\\
\varphi_{uu}&=&0\label{eq:29}\\
2\varphi_{xu}&=&\xi^1_{xx}+\xi^1_{yy}\label{eq:30}\\
\varphi_{yu}&=&\xi^2_{xx}+\xi^2_{yy}\label{eq:31}\\
\varphi_u&=&\xi_x^1=\xi_y^2\label{eq:32}\\
\mu&=&(\xi^1_x-\xi_t^3)f\label{eq:33}\\
\varphi_{tt}&=&f(\varphi_{xx}+\varphi_{yy})\label{eq:34}
\end{eqnarray}
so, we find that
\begin{eqnarray}\nonumber
&&\xi^1(x)=c_1x+c_2y+c_3,\hspace{1cm}\xi^2(t)=c_1y+c_4,\hspace{1cm}\xi^3(t)=a(t),\\
&&\hspace{1cm}\varphi(x,y,u)=c_1u+\beta(x,y),\hspace{1cm}\mu=(c_1-a'(t))f,\label{eq:35}
\end{eqnarray}
with constants $c_1, c_2, c_3$ and $c_4$, also we have
$\beta_{xx}=-\beta_{yy}$.

$\;\;\;\;$We summarize: The class of Eq. (\ref{eq:2}) has an
infinite continuous group of equivalence transformations generated
by the following infinitesimal operators:
\begin{eqnarray}
Y=(c_1x+c_2y+c_3)\frac{\p}{{\p}x}+ (c_1y+c_4)\frac{\p}{{\p}y}+
a(t)\frac{\p}{{\p}t}+(c_1u+\beta(x,y))\frac{\p}{{\p}u}
+(c_1-a'(t))f\frac{\p}{{\p}f}.\label{eq:36}
\end{eqnarray}
Therefore the symmetry algebra of the Burgers' equation
(\ref{eq:2}) is spanned by the vector fields
\begin{eqnarray}
&Y_1=x\frac{\p}{{\p}x}+y\frac{\p}{{\p}y}+t\frac{\p}{{\p}t}+u\frac{\p}{{\p}u}+f\frac{\p}{{\p}f},
\hspace{1cm} Y_2=y\frac{\p}{{\p}x},\hspace{1cm}
Y_3=\frac{\p}{{\p}x},\hspace{1cm}Y_4=\frac{\p}{{\p}y}&\\\label{eq:37}
&Y_5=a(t)\frac{\p}{{\p}t}-a'(t)f\frac{\p}{{\p}f},\hspace{1cm}
Y_{\beta}=\beta(x,y)\frac{\p}{{\p}u}.&\nonumber
\end{eqnarray}

Moreover, in the group of equivalence transformations there are
included also discrete transformations, i.e., reflections
\begin{eqnarray}
t\longrightarrow-t,\hspace{1.5cm}x\longrightarrow-x,\hspace{1.5cm}u\longrightarrow-u,\hspace{1.5cm}
f\longrightarrow-f.\label{eq:38}
\end{eqnarray}
\begin{table}
\caption{Commutation relations satisfied by infinitesimal
generators in (4.38) }\label{table:1}
\vspace{-0.3cm}\begin{eqnarray*}\hspace{-0.75cm}\begin{array}{llllll}
\hline
  [\,,\,]&\hspace{2cm}Y_1 &\hspace{2cm}Y_2  &\hspace{2cm}Y_3  &\hspace{2cm}Y_4  &\hspace{2cm}Y_5 \hspace{2cm}Y_6 \\ \hline
  Y_1    &\hspace{2cm} 0  &\hspace{2cm} 0   &\hspace{1.8cm}-Y_3 &\hspace{1.8cm}-Y_4 &\hspace{2cm}0   \hspace{1.6cm}-Y_6 \\
  Y_2    &\hspace{2cm} 0  &\hspace{2cm} 0   &\hspace{2cm} 0   &\hspace{1.8cm}-Y_3 &\hspace{2cm}0   \hspace{2cm}0\\
  Y_3    &\hspace{2cm}Y_3 &\hspace{2cm} 0   &\hspace{2cm} 0   &\hspace{2cm}0    &\hspace{2cm}0   \hspace{2cm}0\\
  Y_4    &\hspace{2cm}Y_4 &\hspace{2cm} Y_3 &\hspace{2cm} 0   &\hspace{2cm}0    &\hspace{2cm}0   \hspace{2cm}0  \\
  Y_5    &\hspace{2cm} 0  &\hspace{2cm} 0   &\hspace{2cm} 0   &\hspace{2cm}0    &\hspace{2cm}0   \hspace{2cm}0\\
  Y_6    &\hspace{2cm}Y_6 &\hspace{2cm} 0   &\hspace{2cm} 0   &\hspace{2cm}0    &\hspace{2cm}0   \hspace{2cm}0\\
  \hline
\end{array}\end{eqnarray*}
\end{table}
\begin{table}
\caption{Adjoint relations satisfied by infinitesimal generators
in (4.38) }\label{table:1}
\vspace{-0.3cm}\begin{eqnarray*}\hspace{-0.75cm}\begin{array}{llllll}
\hline
[\,,\,]&\hspace{2cm}Y_1&\hspace{2cm}Y_2&\hspace{2cm}Y_3&\hspace{2cm}Y_4&\hspace{2cm}Y_5\hspace{2cm}Y_6
\\ \hline
  Y_1  &\hspace{2cm}Y_1&\hspace{2cm}Y_2&\hspace{2cm}e^sY_3&\hspace{2cm}e^sY_4&\hspace{2cm}Y_5\hspace{2cm}e^sY_6     \\
  Y_2  &\hspace{2cm}Y_1&\hspace{2cm}Y_2&\hspace{2cm}Y_3&\hspace{2cm}Y_4+sY_3&\hspace{2cm}Y_5\hspace{2cm}Y_6    \\
  Y_3  &\hspace{1.5cm}Y_1-sY_3&\hspace{2cm}Y_2&\hspace{2cm}Y_3&\hspace{2cm}Y_4 &\hspace{2cm}Y_5\hspace{2cm}Y_6   \\
  Y_4  &\hspace{1.5cm}Y_1-sY_4&\hspace{1.6cm}Y_2-sY_3&\hspace{2cm}Y_3&\hspace{2cm}Y_4&\hspace{2cm}Y_5\hspace{2cm}Y_6    \\
  Y_5  &\hspace{2cm}Y_1&\hspace{2cm}Y_2&\hspace{2cm}Y_3&\hspace{2cm}Y_4&\hspace{2cm}Y_5\hspace{2cm}Y_6  \\
  Y_6  &\hspace{1.5cm}Y_1-sY_6&\hspace{2cm}Y_2&\hspace{2cm}Y_3&\hspace{2cm}Y_4&\hspace{2cm}Y_5\hspace{2cm}Y_6  \\
  \hline
\end{array}\end{eqnarray*}
\end{table}
\section{Preliminary group classification}
One can observe in many applications of group analysis that most
of extensions of the principal Lie algebra admitted by the
equation under consideration are taken from the equivalence
algebra ${\goth g}_{\EE}$. We call these extensions
$\EE$-extensions of the principal Lie algebra. The classification
of all nonequivalent equations (with respect to a given
equivalence group $G_{\EE}$,) admitting $\EE$-extensions of the
principal Lie algebra is called a preliminary group
classification. Here, $G_{\EE}$ is not necessarily the largest
equivalence group but, it can be any subgroup of the group of all
equivalence transformations.\newline
So, we can take any finite-dimensional subalgebra (desirable as
large as possible) of an infinite-dimensional algebra with basis
(\ref{eq:31}) and use it for a preliminary group classification.
We select the subalgebra ${\goth g}_6$ spanned on the following
operators:
\begin{eqnarray}
&Y_1=x\frac{\p}{{\p}x}+y\frac{\p}{{\p}y}+t\frac{\p}{{\p}t}+u\frac{\p}{{\p}u}+f\frac{\p}{{\p}f},\hspace{1cm}
Y_2=y\frac{\p}{{\p}x},\hspace{1cm}
Y_3=\frac{\p}{{\p}x},\hspace{1cm}
Y_4=\frac{\p}{{\p}y},&\nonumber\\
&Y_5=\frac{\p}{{\p}t}-f\frac{\p}{{\p}f},\hspace{1cm}Y_6=\frac{\p}{{\p}u}.&\label{eq:39}
\end{eqnarray}
The communication relations between these vector fields is given
in Table 1. To each $s$-parameter subgroup there corresponds a
family of group invariant solutions. So, in general, it is quite
impossible to determine all possible group-invariant solutions of
a PDE. In order to minimize this search, it is useful to construct
the optimal system of solutions. It is well known that the problem
of the construction of the optimal system of solutions is
equivalent to that of the construction of the optimal system of
subalgebras \cite{[2],[12]}. Here, we will deal with the
construction of the optimal system of subalgebras of ${\goth
g}_5$.\newline
Let $G$ be a Lie group, with ${\goth g}$ its Lie algebra. Each
element $T\in G$ yields inner automorphism $T_a\longrightarrow
TT_aT^{-1}$ of the group $G$. Every automorphism of the group $G$
induces an automorphism of ${\goth g}$. The set of all these
automorphism is a Lie group called {\it the adjoint group $G^A$}.
The Lie algebra of $G^A$ is the adjoint algebra of ${\goth g}$,
defined as follows. Let us have two infinitesimal generators
$X,Y\in L$. The linear mapping ${\rm
Ad}X(Y):Y\longrightarrow[X,Y]$ is an automorphism of ${\goth g}$,
called {\it the inner derivation of ${\goth g}$}. The set of all
inner derivations ${\rm ad}X(Y)(X,Y\in{\goth g})$ together with
the Lie bracket $[{\rm Ad}X,{\rm Ad}Y]={\rm Ad}[X,Y]$ is a Lie
algebra ${\goth g}^A$ called the {\it adjoint algebra of ${\goth
g}$}. Clearly ${\goth g}^A$ is the Lie algebra of $G^A$. Two
subalgebras in ${\goth g}$ are {\it conjugate} (or {\it similar})
if there is a transformation of $G^A$ which takes one subalgebra
into the other. The collection of pairwise non-conjugate
$s$-dimensional subalgebras is the optimal system of subalgebras
of order $s$. The construction of the one-dimensional optimal
system of subalgebras can be carried out by using a global matrix
of the adjoint transformations as suggested by Ovsiannikov
\cite{[2]}. The latter problem, tends to determine a list (that is
called an {\it optimal system}) of conjugacy inequivalent
subalgebras with the property that any other subalgebra is
equivalent to a unique member of the list under some element of
the adjoint representation i.e. $\overline{{\goth h}}\,{\rm
Ad(g)}\,{\goth h}$ for some ${\rm g}$ of a considered Lie group.
Thus we will deal with the construction of the optimal system of
subalgebras of ${\goth g}_6$.

The adjoint action is given by the Lie series
\begin{eqnarray}
{\rm Ad}(\exp(s\,Y_i))Y_j
=Y_j-s\,[Y_i,Y_j]+\frac{s^2}{2}\,[Y_i,[Y_i,Y_j]]-\cdots,\label{eq:40}
\end{eqnarray}
where $s$ is a parameter and $i,j=1,\cdots,6$. The adjoint
representations of ${\goth g}_6$ is listed in Table 2, it consists
the separate adjoint actions of each element of ${\goth g}_6$ on
all other elements.

{\bf Theorem 4.1.} {\it An optimal system of one-dimensional Lie
subalgebras of general Burgers' equation (\ref{eq:2}) is provided
by those generated by}
\begin{eqnarray*}
&1)&Y^1=Y_1=x{\p}_x+y{\p}_y+t{\p}_t+u{\p}_u+f{\p}_f,\hspace{2.2cm}2)~Y^2=Y_2=y{\p}_x,\\
&3)&Y^3=-Y_4=-{\p}_y,\hspace{5.7cm}4)~Y^4=Y_1+Y_5=x{\p}_x+y{\p}_y+(t+1){\p}_t+u{\p}_u,\\
&5)&Y^5=Y_1-Y_2=(x-y){\p}_x+y{\p}_y+t{\p}_t+u{\p}_u+f{\p}_f,\hspace{0.5cm}6)~Y^6=Y_2-Y_4=y{\p}_x-{\p}_y,\\
&7)&Y^7=-Y_4+Y_6=-{\p}_y+{\p}_u,\hspace{4.1cm}8)~Y^{8}=-Y_4-Y_6=-{\p}_y-{\p}_u,\\
&9)&Y^9=Y_2+Y_5=y{\p}_x+{\p}_t-f{\p}_f,\hspace{3.3cm}10)~Y^{10}=Y_2-Y_5=y{\p}_x-{\p}_t+f{\p}_f,\\
&11)&Y^{11}=Y_2+Y_6=y{\p}_x+{\p}_u,\hspace{4.1cm}12)~Y^{12}=Y_2-Y_6=y{\p}_x-{\p}_u,\\
&13)&Y^{13}=Y_1+Y_2=(x+y){\p}_x+y{\p}_y+t{\p}_t+u{\p}_u+f{\p}_f,\hspace{1mm}14)~Y^{14}=-Y_4+Y_5+Y_6=-{\p}_y+{\p}_t+{\p}_u-f{\p}_f,\\
&15)&Y^{15}=Y_2-Y_4-Y_5+Y_6=y{\p}_x-{\p}_y-{\p}_t+{\p}_u+f{\p}_f,16)~Y^{16}=Y_2-Y_4+Y_6=y{\p}_x-{\p}_y+{\p}_u,\\
&17)&Y^{17}=Y_2-Y_4+Y_5-Y_6=y{\p}_x-{\p}_y+{\p}_t-{\p}_u-f{\p}_f,18)~Y^{18}=Y_2-Y_4-Y_6=y{\p}_x-{\p}_y-{\p}_u,\\
&19)&Y^{19}=Y_1+Y_2+Y_5=(x+y){\p}_x+(t+1){\p}_t+u{\p}_u,\hspace{0.5cm}20)~Y^{20}=Y_2+Y_5+Y_6=y{\p}_x+{\p}_t+{\p}_u-f{\p}_f,\\
&21)&Y^{21}=Y_2+Y_5-Y_6=y{\p}_x+{\p}_t-{\p}_u-f{\p}_f,\hspace{1.6cm}22)~Y^{22}=Y_2-Y_5-Y_6=y{\p}_x-{\p}_t-{\p}_u+f{\p}_f,\\
&23)&Y^{23}=Y_2-Y_5+Y_6=y{\p}_x-{\p}_t+{\p}_u+f{\p}_f,\hspace{1.6cm}24)~Y^{24}=-Y_4-Y_5-Y_6=-{\p}_y-{\p}_t-{\p}_u+f{\p}_f,\\
&25)&Y^{25}=-Y_4-Y_5+Y_6=-{\p}_y-{\p}_t+{\p}_u+f{\p}_f,\hspace{1.2cm}26)~Y^{26}=-Y_4+Y_5-Y_6=-{\p}_y+{\p}_t-{\p}_u-f{\p}_f,\\
&27)&Y^{27}=Y_2-Y_4+Y_5+Y_6=y{\p}_x-{\p}_y+{\p}_t+{\p}_u-f{\p}_f,\\
&28)&Y^{28}=Y_1+Y_2-Y_5=(x+y){\p}_x+y{\p}_y+(t-1){\p}_t+u{\p}_u+2f{\p}_f,\\
&29)&Y^{29}=Y_1-Y_2-Y_5=(x-y){\p}_x+y{\p}_y+(t-1){\p}_t+u{\p}_u+2f{\p}_f,\\
&30)&Y^{31}=Y_1-Y_2+Y_5=(x-y){\p}_x+y{\p}_y+(t+1){\p}_t+u{\p}_u,\\
&31)&Y^{31}=Y_1-Y_5=x{\p}_x+y{\p}_y+(t-1){\p}_t+u{\p}_u+2f{\p}_f\\
&32)&Y^{32}=Y_2-Y_4-Y_5-Y_6=y{\p}_x-{\p}_y-{\p}_t-{\p}_u+f{\p}_f,
\end{eqnarray*}
{\bf Proof.} Let ${\goth g}_6$ is the symmetry algebra of
Eq.~(\ref{eq:2}) with adjoint representation determined in Table 2
and
\begin{eqnarray}
Y=a_1Y_1+a_2Y_2+a_3Y_3+a_4Y_4+a_5Y_5+a_6Y_6,
\end{eqnarray}
is a nonzero vector field of ${\goth g}_6$. We will simplify as
many of the coefficients $a_i;i=1,\ldots,6$, as possible through
proper adjoint applications on $Y$. We follow our aim in the below
easy cases:\newline
{\it Case 1:} \newline
At first, assume that $a_1\neq 0$. Scaling $Y$ if
necessary, we can assume that $a_1=1$ and so we get
\begin{eqnarray}
Y=Y_1+a_2Y_2+a_3Y_3+a_4Y_4+a_5Y_5+a_6Y_6.
\end{eqnarray}
Using the table of adjoint (Table 2) , if we act on $Y$ with ${\rm
Ad}(\exp(a_3Y_3))$, the coefficient of $Y_3$ can be vanished:
\begin{eqnarray}
Y'=Y_1+a_2Y_2+a_4Y_4+a_5Y_5+a_6Y_6.
\end{eqnarray}
Then we apply ${\rm Ad}(\exp(a_4Y_4))$ on $Y'$ to cancel the
coefficient of $Y_4$:
\begin{eqnarray}
Y''=Y_1+a_2Y_2+a_5Y_5+a_6Y_6.
\end{eqnarray}
At last, we apply ${\rm Ad}(\exp(a_6Y_6))$ on $Y''$ to cancel the
coefficient of $Y_6$:
\begin{eqnarray}
Y'''=Y_1+a_2Y_2+a_5Y_5.
\end{eqnarray}
{\it Case 1a:} \\
If $a_2,a_5\neq 0$ then we can make the coefficient of
$Y_2$ and $Y_5$ either $+1$ or $-1$. Thus any one-dimensional
subalgebra generated by $Y$ with $a_2,a_5\neq 0$ is equivalent to
one generated by $Y_1\pm Y_2\pm Y_5$ which introduce parts 19),
28), 29) and  30) of the theorem.\newline
{\it Case 1b:} \newline
 For $a_2=0, a_5\neq0$ we can see that each one-dimensional
subalgebra generated by $Y$ is equivalent to one generated by $
Y_1\pm Y_5$ which introduce parts  4) and  31) of the
theorem.\newline
{\it Case 1c:} \newline
 For $a_2\neq0, a_5=0$, each one-dimensional subalgebra
generated by $Y$ is equivalent to one generated by $Y_1\pm Y_2$
which introduce parts 5) and  13) of the theorem.\newline
{\it Case 1d:} \newline
 For $a_2=0, a_5=0$, each one-dimensional subalgebra
generated by $Y$ is equivalent to one generated by $Y_1$ which
introduce parts  1)  of the theorem.\newline
{\it Case 2:} \newline
 The remaining one-dimensional subalgebras are spanned by
vector fields of the form $Y$ with $a_1=0$. \newline
{\it Case 2a:} \newline
 If $a_4\neq 0$ then by scaling $Y$, we can assume that
$a_4=-1$. Now by the action of ${\rm Ad}(\exp a_3Y_3))$ on $Y$, we
can cancel the coefficient of $Y_3$:
\begin{eqnarray}
\overline{Y}=a_2Y_2-Y_4+a_5Y_5+a_6Y_6.
\end{eqnarray}
Let $a_2\neq0$ then by  scaling $Y$, we can assume that $a_2=1$,
and we have
\begin{eqnarray}
\overline{Y}'=Y_2-Y_4+a_5Y_5+a_6Y_6.
\end{eqnarray}
{\it Case 2a-1:} \newline
Suppose $a_5=a_6=0$, then the one-dimensional subalgebra generated
by $Y$ is equivalent to one generated by $Y_2-Y_4$ which introduce
parts  6). \newline
{\it Case 2a-2:} \newline
Suppose $a_5=0, a_6\neq0$, all of the one-dimensional subalgebra
generated by $Y$ is equivalent to one generated by $Y_2-Y_4\pm
Y_6$ which introduce parts   16) and   18).\newline
{\it Case 2a-3:} \newline
Suppose $a_5\neq0, a_6\neq0$, all of the one-dimensional
subalgebra generated by $Y$ is equivalent to one generated by
$Y_2-Y_4\pm Y_5\pm Y_6$ which introduce parts  15), 17), 27), and
32).

Now if $a_2=0$, we have
\begin{eqnarray}
\overline{Y}''=-Y_4+a_5Y_5+a_6Y_6.
\end{eqnarray}
{\it Case 2a-4:} \newline
Suppose $a_5=a_6=0$, then the one-dimensional subalgebra generated
by $Y$ is equivalent to one generated by $-Y_4$ which introduce
parts 3). \newline
{\it Case 2a-5:} \newline
Suppose $a_5=0, a_6\neq0$, all of the one-dimensional subalgebra
generated by $Y$ is equivalent to one generated by $-Y_4\pm Y_6$
which introduce parts  7) and  8).\newline
{\it Case 2a-6:} \newline
Suppose $a_5\neq0, a_6\neq0$, all of the one-dimensional
subalgebra generated by $Y$ is equivalent to one generated by
$-Y_4\pm Y_5\pm Y_6$ which introduce parts  14), 24), 25) and 26).

{\it Case 2b:} \newline
$~~~~$ Let $a_4=0$ then  $Y$ is in the
form
\begin{eqnarray}
\widehat{Y}=a_2Y_2+a_5Y_5+a_6Y_6.
\end{eqnarray}
Suppose that $a_2\neq 0$ then if necessary we can let it equal to
$1$ and we obtain
\begin{eqnarray}
\widehat{Y}'=Y_2+a_5Y_5+a_6Y_6.
\end{eqnarray}
{\it Case 2b-1:} \newline
Let $a_5=a_6=0$, then $Y_2$ is remained and find 2) section
of the theorem.\newline
{\it Case 2b-2:} \newline
If $a_5\neq0, a_6\neq0$,  then  $\widehat{Y}'$ is equal to
$Y_2\pm Y_5\pm Y_6$. Hence this case suggests part  20), 21), 22)
and 23).\newline
{\it Case 2b-3:} \newline
If $a_5\neq0, a_6=0$,  then  $\widehat{Y}'=Y_2\pm Y_5$ .
Hence this case suggests part 9) and 10).\newline
{\it Case 2b-4:} \newline
If $a_5=0, a_6\neq0$,  then   $Y_2\pm Y_6$ is obtained. So,
this case suggests part 11)  and 12). \newline
There is not any more possible case for studying and the proof is
complete.\hfill\ $\Box$

The coefficients $f$ of Eq. (\ref{eq:2}) depend on the variables
$x,y,u,u_x,u_y$. Therefore, we take their optimal system's
projections on the space $(x,y,u,u_x,u_y,f)$. we have
\begin{eqnarray}\hspace{-0.7cm}
\begin{array}{rlrl}
1)&Z^1=Y^1=x{\p}_x+y{\p}_y+u{\p}_u+f{\p}_f, \hspace{1cm}&17)&Z^{17}=Y^{17}=(x-y){\p}_x+y{\p}_y+u{\p}_u+2f{\p}_f,\\
2)&Z^2=Y^2=y{\p}_x,\hspace{1cm}&18)&Z^{18}=Y^{18}=(x+y){\p}_x+u{\p}_u,\\
3)&Z^3=Y^3=-{\p}_y,\hspace{1cm}&19)&Z^{19}=Y^{19}=y{\p}_x+{\p}_u-f{\p}_f,\\
4)&Z^4=Y^4=(x+y){\p}_x+y{\p}_y+u{\p}_u+f{\p}_f,\hspace{1cm}&20)&Z^{20}=Y^{20}=y{\p}_x-{\p}_u-f{\p}_f,\\
5)&Z^5=Y^5=x{\p}_x+y{\p}_y+u{\p}_u,\hspace{1cm}&21)&Z^{21}=Y^{21}=y{\p}_x-{\p}_u+f{\p}_f,\\
6)&Z^6=Y^6=x{\p}_x+y{\p}_y+u{\p}_u+2f{\p}_f,\hspace{1cm} &22)&Z^{22}=Y^{22}=y{\p}_x+{\p}_u+f{\p}_f,\\
7)&Z^7=Y^7=(x-y){\p}_x+y{\p}_y+u{\p}_u+f{\p}_f,\hspace{1cm} &23)&Z^{23}=Y^{23}=y{\p}_x-{\p}_y+{\p}_u,\\
8)&Z^8=Y^8=y{\p}_x-{\p}_y,\hspace{1cm}&24)&Z^{24}=Y^{24}=y{\p}_x-{\p}_y-{\p}_u,\\
9)&Z^9=Y^9=-{\p}_y+{\p}_u,\hspace{1cm}&25)&Z^{25}=Y^{25}=-{\p}_y+{\p}_u-f{\p}_f,\\
10&Z^{10}=Y^{10}=-{\p}_y-{\p}_u,\hspace{1cm}&26)&Z^{26}=Y^{26}=-{\p}_y-{\p}_u+f{\p}_f,\\
11&Z^{11}=Y^{11}=y{\p}_x-f{\p}_f,\hspace{1cm}&27)&Z^{27}=Y^{27}=-{\p}_y+{\p}_u+f{\p}_f,\\
12)&Z^{12}=Y^{12}=y{\p}_x+f{\p}_f,\hspace{1cm}&28)&Z^{28}=Y^{28}=-{\p}_y-{\p}_u-f{\p}_f,\\
\end{array}
\end{eqnarray}
\begin{eqnarray}\hspace{-0.7cm}
\begin{array}{rlrl}
13)&Z^{13}=Y^{13}=y{\p}_x+{\p}_u,\hspace{1cm}&29)&Z^{29}=Y^{29}=y{\p}_x-{\p}_y+{\p}_u-f{\p}_f,\\
14)&Z^{14}=Y^{14}=y{\p}_x-{\p}_u,\hspace{1cm}&30)&Z^{30}=Y^{30}=y{\p}_x-{\p}_y-{\p}_u+f{\p}_f,\\
15)&Z^{15}=Y^{15}=(x-y){\p}_x+y{\p}_y+u{\p}_u,\hspace{1cm}&31)&Z^{31}=Y^{31}=y{\p}_x-{\p}_y+{\p}_u+f{\p}_f,\\
16)&Z^{16}=Y^{16}=(x+y){\p}_x+y{\p}_y+u{\p}_u+2f{\p}_f,\hspace{1cm}&32)&Z^{32}=Y^{32}=y{\p}_x-{\p}_y-{\p}_u-f{\p}_f,
\end{array}
\end{eqnarray}
%
{\bf Proposition 4.2.} {\it Let ${\goth g}_m:=\langle Y_1, \ldots,
Y_m\rangle$, be an $m$-dimensional algebra. Denote by $Y^i (i=1,
\ldots, r, 0<r\leq m, r\in{\Bbb N})$ an optimal system of
one-dimensional subalgebras of ${\goth g}_m$ and by $Z^i\, (i =
1,\cdots, t, 0<t\leq r, t\in{\Bbb N})$ the projections of $Y^i$,
i.e., $Z^i = {\rm pr}(Y^i)$. If equations
\begin{eqnarray}
f = \Phi(x,y,u,u_x,u_y),\label{eq:18}
\end{eqnarray}
are invariant with respect to the optimal system $Z^i$ then the
equation
\begin{eqnarray}
u_t = \Phi(x,y,u,u_x,u_y)(u_{xx}+u_{yy}),\label{eq:19}
\end{eqnarray}
admits the operators $X^i=$ projection of $Y^i$ on $(t,x,y,u,u_x,u_y)$.}

{\bf Proposition 4.3.} {\it Let Eq. (\ref{eq:19}) and the equation
\begin{eqnarray}
u_t =  \Phi'(x,y,u,u_x,u_y)(u_{xx}+u_{yy}),\label{eq:20}
\end{eqnarray}
be constructed according to Proposition 4.2. via optimal systems
$Z^i$ and ${Z^i}'$, respectively. If the subalgebras spanned on
the optimal systems $Z^i$ and ${Z^i}'$, respectively, are similar
in ${\goth g}_m$, then Eqs. (\ref{eq:19}) and (\ref{eq:20}) are
equivalent with respect to the equivalence group $G_m$, generated
by ${\goth g}_m$. }

Now we apply Proposition 4.2. and  Proposition 4.3. to the optimal
system (\ref{eq:17}) and obtain all nonequivalent Eq. (\ref{eq:2})
admitting $\EE$-extensions of the principal Lie algebra ${\goth
g}$, by one dimension, i.e., equations of the form (\ref{eq:2})
such that they admit, together with the one basic operators
(\ref{eq:21}) of ${\goth g}$, also a second operator $X^{(2)}$.
For every case, when this extension occurs, we indicate the
corresponding coefficients $f, g$ and the additional operator
$X^{(2)}$.

 We perform the algorithm passing from operators
$Z^i\,(i=1,\cdots,32)$ to $f$ and $X^{(2)}$ via the following
example. \newline
Let consider the vector field
\begin{eqnarray}
Z^{32}=y{\p}_x-{\p}_y-{\p}_u-f{\p}_f,\label{eq:21}
\end{eqnarray}
then the characteristic equation corresponding to $Z^6$ is
\begin{eqnarray}
{dx\over y}={dy\over-1}=\frac{du}{-1}=\frac{df}{-f},
\end{eqnarray}
and can be taken in the form
\begin{eqnarray}
I_1=u+{x\over y},\hspace{5mm}I_2=e^{x\over y}f.
\end{eqnarray}
From the invariance equations we can write
\begin{eqnarray}
I_2=\Phi(I_1),
\end{eqnarray}
it follows that
\begin{eqnarray}
f=e^{-{x\over y}}\Phi(\lambda),
\end{eqnarray}
where $\lambda=I_1$.

From Proposition 4.2. applied to the operator $Z^6$ we obtain the
additional operator $X^{(2)}$
\begin{eqnarray}
y{\p}_x-{\p}_y+{\p}_t-{\p}_u.
\end{eqnarray}
After similar calculations applied to all operators (\ref{eq:17})
we obtain the following result (Table 3) for the preliminary group
classification of Eq. (\ref{eq:2}) admitting an extension ${\goth
g}_3$ of the principal Lie algebra ${\goth g}_1$.
\section{Conclusion}
In this paper, following the classical Lie method, the preliminary
group classification for the class of heat equation (\ref{eq:2})
and investigated the algebraic structure of the symmetry groups
for this equation, is obtained. The classification is obtained by
constructing an optimal system with the aid of Propositions 4.2.
and 4.3.. The result of the work is summarized in Table 3. Of
course it is also possible to obtain the corresponding reduced
equations for all the cases in the classification reported in
Table 3.

\begin{table}
\centering{\caption{The result of the
classification}}\label{table:3} \vspace{-0.35cm}\begin{eqnarray*}
\hspace{-0.75cm}\begin{array}{l l l l l l} \hline
  N       &\hspace{1cm} Z     &\hspace{1.1cm} \mbox{Invariant}  &\hspace{1cm} \mbox{Equation}
  &\hspace{1cm} \mbox{Additional operator}\,X^{(2)} \\ \hline
  1       &\hspace{1cm} Z^1   &\hspace{1.1cm} {u\over x}  &\hspace{1cm}u_t=x\Phi(u_{xx}+u_{yy})
  &\hspace{1cm} x{\p}_x+y{\p}_y+t{\p}_t+u{\p}_u \\
  2       &\hspace{1cm} Z^2   &\hspace{1.1cm} u   &\hspace{1cm}u_t=\Phi(u_{xx}+u_{yy})
          &\hspace{1cm} y{\p}_x \\
  3       &\hspace{1cm} Z^3   &\hspace{1.1cm} u  &\hspace{1cm}u_t=\Phi(u_{xx}+u_{yy})
  &\hspace{1cm} -{\p}_y\\
  4        &\hspace{1cm} Z^4    &\hspace{1.1cm} {u\over x+y}  &\hspace{1cm}u_t=y\Phi(u_{xx}+u_{yy})
  &\hspace{1cm} (x+y){\p}_x+y{\p}_y+u{\p}_u \\
  5        &\hspace{1cm} Z^5    &\hspace{1.1cm} {u\over x}&\hspace{1cm}u_t=\Phi(u_{xx}+u_{yy})
          &\hspace{1cm} x{\p}_x+y{\p}_y+(t+1){\p}_t+u{\p}_u \\
  6        &\hspace{1cm} Z^6   &\hspace{1.1cm} {u\over x}&\hspace{1cm}u_t=x^2\Phi(u_{xx}+u_{yy})
  &\hspace{1cm} x{\p}_x+y{\p}_y+(t-1){\p}_t+u{\p}_u \\
  7        &\hspace{1cm} Z^7   &\hspace{1.1cm} {u\over x-y}&\hspace{1cm}u_t=(x-y)\Phi(u_{xx}+u_{yy})
  &\hspace{1cm}(x-y){\p}_x+y{\p}_y+t{\p}_t+u{\p}_u \\
  8        &\hspace{1cm} Z^8   &\hspace{1.1cm} u &\hspace{1cm}u_t=\Phi(u_{xx}+u_{yy})
  &\hspace{1cm} y{\p}_x-{\p}_y \\
  9        &\hspace{1cm} Z^9   &\hspace{1.1cm}  u&\hspace{1cm}u_t=\Phi(u_{xx}+u_{yy})
  &\hspace{1cm} -{\p}_y+{\p}_u \\
  10        &\hspace{1cm} Z^{10}   &\hspace{1.1cm} x &\hspace{1cm}u_t=\Phi(u_{xx}+u_{yy})
  &\hspace{1cm} -{\p}_y-{\p}_u \\
  11       &\hspace{1cm} Z^{11}   &\hspace{1.1cm} u &\hspace{1cm}u_t=e^{-{x\over y}}\Phi(u_{xx}+u_{yy})
  &\hspace{1cm} y{\p}_x+{\p}_t \\
  12       &\hspace{1cm} Z^{12}   &\hspace{1.1cm} u &\hspace{1cm}u_t=e^{x\over y}\Phi(u_{xx}+u_{yy})
  &\hspace{1cm} y{\p}_x-{\p}_t \\
  13       &\hspace{1cm} Z^{13}   &\hspace{1.1cm}  u-{x\over y} &\hspace{1cm}u_t=\Phi(u_{xx}+u_{yy})
  &\hspace{1cm} y{\p}_x+{\p}_u \\
  14       &\hspace{1cm} Z^{14}   &\hspace{1.1cm} u+{x\over y} &\hspace{1cm}u_t=\Phi(u_{xx}+u_{yy})
  &\hspace{1cm} y{\p}_x-{\p}_u \\
  15       &\hspace{1cm} Z^{15}   &\hspace{1.1cm} {u\over x-y}&\hspace{1cm}u_t=\Phi(u_{xx}+u_{yy})
  &\hspace{1cm} (x-y){\p}_x+y{\p}_y+(t+1){\p}_t+u{\p}_u \\
  16       &\hspace{1cm} Z^{16}   &\hspace{1.1cm} {u\over x+y}&\hspace{1cm}u_t=(x+y)^2\Phi(u_{xx}+u_{yy})
  &\hspace{1cm} (x+y){\p}_x+y{\p}_y+(t-1){\p}_t+u{\p}_u \\
  17       &\hspace{1cm} Z^{17}   &\hspace{1.1cm} {u\over x-y}&\hspace{1cm}u_t=(x-y)^2\Phi(u_{xx}+u_{yy})
  &\hspace{1cm} (x-y){\p}_x+y{\p}_y+(t-1){\p}_t+u{\p}_u \\
  18       &\hspace{1cm} Z^{18}   &\hspace{1.1cm} {u\over x+y}&\hspace{1cm}u_t=\Phi(u_{xx}+u_{yy})
  &\hspace{1cm} (x+y){\p}_x+(t+1){\p}_t+u{\p}_u \\
  19       &\hspace{1cm} Z^{19}   &\hspace{1.1cm}u-{x\over y}&\hspace{1cm}u_t=e^{-{x\over y}}\Phi(u_{xx}+u_{yy})
  &\hspace{1cm} y{\p}_x+{\p}_t+{\p}_u\\
  20       &\hspace{1cm} Z^{20}   &\hspace{1.1cm} u+{x\over y}&\hspace{1cm}u_t=e^{-{x\over y}}\Phi(u_{xx}+u_{yy})
  &\hspace{1cm} y{\p}_x+{\p}_t-{\p}_u \\
  21       &\hspace{1cm} Z^{21}   &\hspace{1.1cm} u+{x\over y}&\hspace{1cm}u_t=e^{x\over y}\Phi(u_{xx}+u_{yy})
  &\hspace{1cm} y{\p}_x-{\p}_t-{\p}_u\\
  22       &\hspace{1cm} Z^{22}   &\hspace{1.1cm} u-{x\over y}&\hspace{1cm}u_t=e^{x\over y}\Phi(u_{xx}+u_{yy})
  &\hspace{1cm} y{\p}_x-{\p}_t+{\p}_u\\
  23       &\hspace{1cm} Z^{23}   &\hspace{1.1cm} u-{x\over y}&\hspace{1cm}u_t=\Phi(u_{xx}+u_{yy})
  &\hspace{1cm} y{\p}_x-{\p}_y+{\p}_u\\
  24       &\hspace{1cm} Z^{24}   &\hspace{1.1cm}u+{x\over y}&\hspace{1cm}u_t=e^y\Phi(u_{xx}+u_{yy})
  &\hspace{1cm} y{\p}_x-{\p}_y-{\p}_u \\
  25       &\hspace{1cm} Z^{25}   &\hspace{1.1cm} u+y&\hspace{1cm}u_t=\Phi(u_{xx}+u_{yy})
  &\hspace{1cm} -{\p}_y+{\p}_t+{\p}_u \\
  26       &\hspace{1cm} Z^{26}   &\hspace{1.1cm} u-y&\hspace{1cm}u_t=e^{-y}\Phi(u_{xx}+u_{yy})
  &\hspace{1cm} -{\p}_y-{\p}_t-{\p}_u \\
  27       &\hspace{1cm} Z^{27}   &\hspace{1.1cm} u+y&\hspace{1cm}u_t=e^{-y}\Phi(u_{xx}+u_{yy})
  &\hspace{1cm} -{\p}_y-{\p}_t+{\p}_u \\
  28       &\hspace{1cm} Z^{28}   &\hspace{1.1cm} u-y &\hspace{1cm}u_t=e^y\Phi(u_{xx}+u_{yy})
  &\hspace{1cm} -{\p}_y+{\p}_t-{\p}_u \\
  29       &\hspace{1cm} Z^{29}   &\hspace{1.1cm} u+{x\over y}&\hspace{1cm}u_t=e^{-{x\over y}}\Phi(u_{xx}+u_{yy})
  &\hspace{1cm} y{\p}_x-{\p}_y+{\p}_t+{\p}_u \\
  30       &\hspace{1cm} Z^{30}   &\hspace{1.1cm} u+{x\over y}&\hspace{1cm}u_t=e^{x\over y}\Phi(u_{xx}+u_{yy})
  &\hspace{1cm} y{\p}_x-{\p}_y-{\p}_t-{\p}_u \\
  31       &\hspace{1cm} Z^{31}   &\hspace{1.1cm} u-{x\over y}&\hspace{1cm}u_t=e^{x\over y}\Phi(u_{xx}+u_{yy})
  &\hspace{1cm} y{\p}_x-{\p}_y-{\p}_t+{\p}_u \\
  32       &\hspace{1cm} Z^{32}   &\hspace{1.1cm} u+{x\over y}&\hspace{1cm}u_t=e^{-{x\over y}}\Phi(u_{xx}+u_{yy})
  &\hspace{1cm} y{\p}_x-{\p}_y+{\p}_t-{\p}_u \\
  \hline
\end{array}\end{eqnarray*}
\end{table}
\vspace{17cm}

\end{document}